\documentclass[10pt,twoside]{article}
\usepackage{Latex-document}
\usepackage{amsmath,amssymb,amsthm,amscd}
\numberwithin{equation}{section} 

\def\p{\partial}
\def\b{\bar}

\def\o{\omega}

\def\RR{{\mathbb R}}

\newtheorem{prop}{Proposition}[section]
\newtheorem{theo}[prop]{Theorem}

\newtheorem{cor}[prop]{Corollary}

\newtheorem{defi}[prop]{Definition}

\newtheorem{q}[prop]{Question}
\def\begeq{\begin{equation}}
\def\endeq{\end{equation}}
\def\and{\quad{\rm and}\quad}

\let\lra=\longrightarrow

\def\mapright\#1{\,\smash{\mathop{\lra}\limits^{\#1}}\,}

\renewcommand{\thesection}{\arabic{section}}
\renewcommand{\thesubsection}{\thesection.\arabic{subsection}.}

\title{\bf Recent Progress in K\"ahler Geometry\thanks{Paritally supported by NSF research grant DMS-0110321
(2001-2004).}\vskip 6mm}
\author{Xiuxiong Chen\thanks{Department of Mathematics, Princeton University.
Department of Mathematics, University of Wisconsin at Madison,
USA. E-mail: xiu@math.princeton.edu}\vspace*{-0.5cm}}
\date{\vspace{-8mm}}

\begin {document}
 \maketitle

\thispagestyle{first} \setcounter{page}{273}

\begin{abstract}

\vskip 3mm

In recent years, there are many progress made in K\"ahler geometry. In particular, the  topics related to the
problems of the existence and uniqueness of extremal K\"ahler metrics,  as well as obstructions to the existence
of such metrics  in general K\"ahler manifold.  In this talk, we will report some recent developments in this
direction.  In particular, we will discuss the progress recently obtained in understanding the metric structure of
the infinite dimensional space of Kaehler potentials, and their applications to the problems mentioned above. We
also will discuss some recent on Kaehler Ricci flow.

\vskip 4.5mm

\noindent {\bf 2000 Mathematics Subject Classification:} 53, 35.

\noindent {\bf Keywords and Phrases:} Extremal K\"ahler metrics, K\"ahler-Einstein metrics, Holomorphic vector
field, Holomorphic invariant, K\"ahler Ricci flow.
\end{abstract}

\vskip 12mm

 In the last few years, we have witnessed a rapid progress in
 K\"ahler geometry. In particular, the topic related to the
 existence, to the uniqueness of extremal K\"ahler metrics,  and to obstructions to the existence of such metrics.
  In this talk, we will give a brief survey of these exciting progress made in this direction.

\vspace*{-0.5mm}

\subsection{Some background}

\vskip-5mm \hspace{5mm}

Let $(M,\omega)$ be a polarized $n$-dimensional compact K\"ahler
manifold, where $\omega$ is a K\"ahler form  on $M$. In local
coordinates $z_1, \cdots, z_n$, we have
\[
\omega = \sqrt{-1} \displaystyle \sum_{i,j=1}^n\;g_{i
\overline{j}} d\,z^i\wedge d\,z^{\overline{j}}  > 0,
\]
where $\{g_{i\overline {j}}\}$ is a positive definite Hermitian
matrix function. The K\"ahler condition requires that $\omega$ is
a closed positive (1,1)-form. The K\"ahler metric corresponding to
$\omega$ is given by
\[
 g_\omega = \displaystyle \sum_{\alpha, \beta = 1}^n \; {g}_{\alpha \overline{\beta}} \;
d\,z^{\alpha}\;\otimes d\, z^{ \overline{\beta}}.
\]
For simplicity, in the following, we will often denote by $\omega$
the corresponding K\"ahler metric. The K\"ahler class of $\omega$
is its cohomology class $[\omega]$ in $H^2(M,\RR).\;$ It follows
from the Hodge-Dolbeault theorem that any other K\"ahler metric in
the same K\"ahler class is of the form
\[
\omega_{\varphi} = \omega + \sqrt{-1} \displaystyle
\sum_{i,j=1}^n\; {{\partial^2 \varphi}\over {\partial z^i \partial
z^{\overline{j}}}}
> 0
\]
for some real valued function $\varphi$ on $M.\;$

Given a K\"ahler metric $\omega$, its volume form  is
\[
 {1 \over n!} \omega^n = \left(\sqrt{-1} \right)^n \det\left(g_{i \overline{j}}\right)
 d\,z^1 \wedge d\,z^{\overline{1}}\wedge \cdots \wedge d\,z^n \wedge d\,z^{\overline{n}}.
\]
Its Ricci (curvature) form is: \[ Ric(\omega) = \sqrt{-1} R_{i \b
j} \; d\,w_i \;d\,\b w_j  = - \sqrt{-1} \p \b \p \;\log\; \det\;
\omega^n.
\]
Note also that $R(\omega) = g^{i\b j} R_{i \b j}$ corresponds to
one half times the scalar curvature as it is usually defines in
Riemannian geometry.  We say that the first Chern class of $M$ is
positive of negative definite, if there exists a real valued
function $\psi$ on $M$ such that $ R_{i \b j} + {{\p^2
\psi}\over{\p w_i \p \b w_j}}$ is, respectively, positive of
negative definite. A K\"ahler metric is K\"ahler-Einstein, if the
Ricci form is proportional to the K\"ahler form by a constant
factor.  A K\"ahler metric is called extremal in the sense of E.
Calabi \cite{calabi82}, if it is a critical point of the
functional $\displaystyle \int_M\; |Ric(\omega)|^2 \; \omega^n,\;$
or, equivalently, if the complex gradient vector field of the
scalar curvature function $
  g^{\alpha \b \beta} (\omega) {{\p R(\omega)}\over {\p \b
  z_\beta}} {{\p}\over {\p z_\alpha}} \;$
is a holomorphic vector field.

\vspace*{-0.5mm}

\subsection{Existence of extremal K\"ahler metrics}

\vskip-5mm \hspace{5mm}

It is well known that a K\"ahler-Einstein metric
satisfies a Monge-Ampere equation
\[
   \log \det {{\omega_\varphi^n}\over {\omega^n}} = -\lambda\;
   \varphi + h_\omega
\]
where $[Ric(\omega)] = \lambda \;[\omega]$ and
\[
  Ric(\omega) - \lambda \;\omega = i\p \b \p\; h_\omega.
\]
 In Calabi's work in the 1950s, he  made conjectures
about the existence of K\"ahler-Einstein metrics on compact
K\"ahler manifolds with
 definite first  Chern class. In 1976, Aubin and Yau independently obtained
 existence when the first Chern class is negative. Around the same time, Yau proved
 also the existence of a K\"ahler-Einstein metric when the first Chern class vanishes.
 This is a celebrated work;   and any K\"ahler manifold admit such a metric is called ``Calabi-Yau" manifold.
    The positive case remains
 open, but significant progress has been made in the last two decades.
  G. Tian proved in \cite{Tian90} the existence
 of K\"ahler-Einstein metrics on any complex surface with positive first
 Chern class and reductive automorphism group.  In 1997, Tian \cite{Tian97} proved that existence of
K\"ahler-Einstein metrics with positive scalar curvature is
equivalent to an analytic stability.  It remains open how this
analytic stability follows from certain algebraic stability in
geometric invariant theory.\\

 The construction of complete non-compact Calabi-Yau manifolds has also enjoyed
 a good deal of success through the work of Calabi, Tian and Yau,
 Anderson, Kronheimer, LeBrun, Joyce and many others. These non-compact
 metrics are related to manifolds with $G_2$ and $Spin(7)$ holonomy, which are
 important in M-theory.\\

   A lot of effort has also gone into constructing special or
 explicit examples of K\"ahler-Einstein metrics and extremal K\"ahler metrics.
 The same is true for hyperkaehler metrics as well. Counter examples to the
 existence of extremal metrics have given by Levine, Burns-De Bartolomeis,
 and LeBrun.\\

There has not been much  progress made on  the existence of
extremal metrics in general. One of the possible reasons is the
lack of maximum principle for non-linear equations of 4th order. A
general existence result, even in complex surfaces, will be highly
interesting.

\vspace*{-0.5mm}

\subsection{Obstructions}

\vskip-5mm \hspace{5mm}

      In 1983, A. Futaki \cite{futaki83}
 introduced a complex character ${\cal F}(X,[\omega])$ on the complex Lie algebra of all holomorphic
 vector fields $X$ in $M,\;$ depending only on the K\"ahler class $[\omega],$ and show that
 its vanishing is a necessary condition
 for the existence of a K\"ahler-Einstein metric on the manifold.
In 1985, E. Calabi\cite{calabi85} generalized  Futaki's result to
cover the more general case of any extremal K\"ahler metric: the
generalized Futaki invariant of a given K\"ahler class is zero or
not, according to whether any extremal metric in that class has
constant scalar curvature or not. S. Bando  also obtained some
generalizations of the Futaki  invariant. More recently, a finite
family of obstructions was introduced
 in \cite{chentian001}.
For any holomorphic vector field $X$ inducing the trivial
translation on the Albanese torus there exists a complex valued
potential function $\theta_{X,\omega},$ uniquely determined up to
additive constants, defined by the equation: $
   L_{X} \omega = \sqrt{-1} \p \bar \p \theta_X $(
Here $L_X$ denote the Lie derivative along vector field $X.$).
Now, for each $k=0,1,\cdots n,$ define the functional
$\Im_k(X,\omega)$ by\footnote{This is a formula for canonical
K\"ahler class. For general K\"ahler class, see
\cite{chentian001}.}
\[
\begin{array}{ll}
\Im_k (X,\omega)  = (n-k) \displaystyle \int_M\theta_X\;
{\omega}^n\\
\qquad  + \displaystyle \int_M \left( (k+1) \Delta \theta_X
\;{{\rm Ric}(\omega)}^{k}\wedge {\omega}^{n-k} - (n-k)\; \theta_X
\;{{\rm Ric}(\omega)}^{k+1} \wedge {\omega}^{n-k-1}\right ).
\end{array}
\]
Here and elsewhere,  $\Delta_\omega$ denotes the one half times
the Laplacian-Beltrami operator of the induced Riemannian
structure $\omega.\;$

The next theorem assures that the above integral gives rise to a
holomorphic invariant.

\begin{theo} {\rm \cite{chentian001}}The integral $\Im_k (X,\omega)$ is independent of choices of
K\"ahler metrics in the K\"ahler class $[\omega]$, that is, $\Im_k
(X,\omega)=\Im_k (X,\omega')$ so long as the K\"ahler forms
$\omega$ and $\omega'$ represent the same K\"ahler class. Hence,
the integral $\Im_k (X,\omega)$ is a holomorphic invariant, which
will be denoted by $\Im_k (X,[\omega])$. Note that $\Im_0$ is the
usual Futaki invariant.
\end{theo}

\vspace*{-0.5mm}

\subsection{Uniqueness of extremal K\"ahler metrics}

\vskip-5mm \hspace{5mm}

     We now turn to the uniqueness of extremal metrics. In the 1950s, Calabi
 used the maximum principle to prove the uniqueness of K\"ahler-Einstein metrics
 when the first Chern class is non-positive. In 1987, Mabuchi introduced the
``K-energy", which is essentially a potential function for the
constant
 scalar curvature metric equation. Using the K-energy, he and
 Bando \cite{Bando87} proved that the uniqueness of K\"ahler-Einstein metric up to
 holomorphic transformations when the first Chern
 class is positive.  Recently, Tian and X. H.
Zhu proved that the uniqueness of K\"{a}hler-Ricci Soliton on any
K\"{a}hler manifolds with positive first Chern class.
\begin{theo} {\rm \cite{TianZhu98}, \cite{TianZhu00}} The K\"ahler Ricci soliton of a K\"ahler manifold $M$
is unique modulo the automorphism subgroup $Aut_r(M);$ more
precisely, if $\omega_1, \omega_2$ are two K\"ahler Ricci solitons
with respect to a holomorphic vector field X, i.e., they satisfies
\begin{equation} Ric(\omega_i) - \omega_i = {\cal L}_X (\omega_i),
\qquad {\rm where \;} i = 1,2.
\end{equation}
Then there are automorphism $\sigma$ in $Aut^o(M)$ and $\tau$ in
$Aut_r(M)$ such that $\sigma_*^{-1} X \in \eta_r (M) $ and
$\sigma^* \omega_2 = \tau^* \sigma^* \omega_1,$ where $\eta_r(M)$
denotes the Lie algebra of $Aut_r(M). $ In fact, $\sigma_*^{-1}X$
lies in the center of $\eta_r(M).\;$ Moreover, this vector field
$X$ is unique up to conjugations.
\end{theo}

Following a program of Donaldson (which will be explained in
Subsection 0.7), we proved in 1998 \cite{chen991} that the
uniqueness for constant scalar curvature metric in any K\"ahler
class when $C_1 < 0$ along with some other interesting results:

\begin{theo} {\rm \cite{chen991}}
If the first Chern class is strictly negative, then the extremal
K\"{a}hler metric is unique in each K\"{a}hler class. Moreover,
the K energy must have a uniform lower bound if there exists an
extremal K\"ahler metric in that K\"ahler class.
\end{theo}

Very recently, Donaldson proved a beautiful theorem which states
\begin{theo}{\rm \cite{Dona02}} For algebraic K\"ahler class with no non-trivial
holomorphic vector field, the  constant scalar curvature metric is
unique.
\end{theo}
The two theorems overlaps in a lot cases, but mutually
non-inclusive.

\vspace*{-0.5mm}

\subsection{Lower bound of the K energy}

\vskip-5mm \hspace{5mm}

 According to T. Mabuchi and S.
Bando\cite{Bando87},
 the existence of a lower bound of the K energy is
 a necessary condition for the existence of
K\"{a}hler-Einstein metrics in the first Chern class. Tian
\cite{Tian97}  showed that in a K\"{a}hler manifold with positive
first Chern class and no non-trivial holomorphic fields, the
K\"{a}hler-Einstein metric exists if and only if the Mabuchi
functional is proper.  When the first Chern class is negative,
making use of Tian's explicit formulation \cite{Tian97}, a simple
idea in \cite{chen993} reduces a lower bound of the K energy to
the existence of critical point for the following convex
functional:

\[
  J(\varphi) = - \displaystyle \sum_{p=0}^{n-1}  {1 \over { (p+1)! (n-p-1)!}} \displaystyle  \int_V \; \varphi\; Ricci(\omega_0) \wedge \omega_0^{n-p-1} \; (\partial \overline{\partial} \varphi
  )^p,
\]
where $Ricci(\omega_0) < 0.\;$  In complex surfaces, we solves
this existence problem completely, which leads to the following
 interesting result:

\begin{theo}{\rm  \cite{chen993}}  Suppose $ {\rm dim} V = 2 $ and $C_1 (V) < 0.\;$ For any K\"ahler class $[\omega_0]$,
if  $ 2\;{{[-C_1(V)] \cdot [\omega_0]} \over
{[\omega_0]^2}}\;\;[\omega_0] + [C_1(V)]  > 0,\; $  then the K
energy has a lower bound in this K\"ahler class.
\end{theo}

It will be very interesting to generalize this result to higher
dimensional K\"ahler manifold.

\vspace*{-0.5mm}

\subsection{Donaldson's program }

\vskip-5mm \hspace{5mm}

     Mabuchi defined in \cite{Ma87} a Weil-Petersson type metric on the space of
 K\"ahler potentials in a fixed K\"ahler class.
 Consider the space of K\"ahler potentials
\[
{\cal H} = \{ \varphi \mid \omega_{\varphi} = \omega + \b \p \p
\varphi > 0,\;{\rm on} \; M\}.
 \]

 A tangent vector in $\cal {H}$ is
just a real valued function in $M.\;$ For any  vector $\psi \in
T_{\varphi} \cal {H}, $ we define the length of this vector as:
\[
\|\psi\|^2_{\varphi} =\int_{V}\psi^2\;d\;\mu_{\varphi}.
\]
It is easy to see that the geodesic equation for this metric is
\[
  \varphi''(t) - g^{\alpha\b \beta}_{\varphi} {{\p \varphi'}\over
  {\p w_\alpha}} {{\p \varphi'}\over
  {\p w_{\b \beta}}} \; = \;0,
\]
where ${g_{\alpha\b \beta}}_{\varphi} = {g_0}_{\alpha \b \beta} +
{{\p^2 \varphi}\over {\p w_{\alpha} \p w_{\b \beta}}} > 0.\; $ It
is first observed (cf.  Semmes S. \cite{Semmes92} )that one can
complexified the $t$ variable, denoted it by  $w_{n+1}.\;$ Then,
the geodesic equation becomes a homogenous complex Monge-Ampere
equation:

\begin{equation}
  \det\; \left( g_{0,i\overline{j}} + {{\partial^2 \varphi}\over{\partial w_{i}
\overline{\partial w_j}}} \right)_{(n+1)(n+1)} \;\; =\;\; 0,
\qquad {\rm on} \; \Sigma \times M. \label{eq:monge0}
\end{equation}
Here $\Sigma = [0,1] \times S^1.\;$  It turns out that we don't
need to restrict to this special case. For any Riemann surface
$\Sigma$ with boundary, and for any $C^\infty$ map $\varphi_0$
from $\p \Sigma$ to $\cal H, $ one can always ask the following
existence problem: \begin{q} {\rm (Donaldson\cite{Dona96})}\rm
For any smooth map $\varphi_0: \p \Sigma \rightarrow \cal H,$ does
there exists a smooth map $\varphi: \Sigma \rightarrow \cal H$
which satisfies the Homogenous Monge-Ampere equation
\ref{eq:monge0} such that $\varphi = \varphi_0$ in $\p \Sigma$ ?
\end{q}
\begin{theo}{\rm  (X. Chen \cite{chen991})} For any smooth map
$\varphi_0: \p \Sigma \rightarrow \cal H,$ there always exists a
$C^{1,1}$ map $\varphi: \Sigma \rightarrow \overline{\cal H} $
which solves the Homogenous Monge-Ampere equation \ref{eq:monge0}
such that $\varphi = \varphi_0$ in $\p \Sigma$.
\end{theo}

 An important conjecture by
 Donaldson in  \cite{Dona96} was that the space of K\"ahler potentials is a metric
 space which is path-connected with respect to this
 Weil-Petersson metric. This conjecture was complete verified
 here.

\begin{theo}{\rm \cite{chen991}}
The space $\cal H$ is a genuine metric space:  the minimal
distance between any two K\"{a}hler metrics is realized by the
unique $C^{1,1}$ geodesic; and the length of this geodesic is
positive.

\end{theo}
Collaborating with E. Calabi, we proved the following
\begin{theo}{\rm \cite{chen992}}
  The space  $\cal H$
 in a fixed K\"{a}hler class is a non-positively curved  space in the sense of Alexandrov: Suppose A, B, C are three
smooth points in $\cal H$  and $P_{\lambda}$
 is a geodesic interpolation point for $ 0 \leq \lambda \leq 1$:
 the distance from $P$ to $B$ and $C$ are respectively $\lambda
 d(B,C)$ and $ (1-\lambda) d(B,c) $ \footnote{In affine space, this means $
P_{\lambda} = \lambda\; B + (1 - \lambda ) \;C.$}.  Then the
following inequality holds:
\[
d(A, P_{\lambda})^2 \leq (1-\lambda) d(A, B)^2 + \lambda d (A,
C)^2 - \lambda \cdot (1-\lambda) d(B,C)^2.
\]
\end{theo}

  \begin{theo}{\rm  \cite{chen992}}  Given
 any two K\"ahler potentials  $\varphi_1$ and
$\varphi_2$ in $\;\cal H$ and a smooth curve $\varphi(t), 0 \leq t
\leq 1 $ which connects them  in $\cal H.\;$ Suppose $\varphi(s,t)
$ are the family of curves under the Calabi flow   and suppose
that $L(s)$ is the length of this curve at time $s.\;$ Then
\[
  {{d \, L}\over {d\, s}} = -  \int_0^1 \left(\int_M |D
   {{\partial \varphi} \over {\partial t}}|_{\varphi(s,t)}^2 \; \omega_{\varphi(s,t)}^n
    \cdot \sqrt{ \int_M |{{\partial \varphi} \over {\partial t}}|^2\; \; \omega_{\varphi(s,t)}^n}^{-{1\over 2}} \right)\; d\,t, \]
where $D$ is the 2nd order Lichernowicz operator. For any smooth
function $f$ in $V,\; D(f) = \displaystyle \sum_{\alpha,
\beta=1}^n\;f_{,\alpha \beta} dz^{\alpha} \otimes d z^{\beta}\;$
where $f_{,\alpha \beta}$ is the second covariant derivatives of
$f.\;$
\end{theo}

\vspace*{-0.5mm}

\subsection{The Calabi flow and the K\"ahler Ricci flow}

\vskip-5mm \hspace{5mm}

In a sequence of papers \cite{chen994}, and  \cite{chentian001}
\cite{chentian002}, we develop some new techniques in attacking
the convergence problems for the geometric flow, in particular,
the Calabi flow and the K\"ahler Ricci flow. The main ideas are to
find a set of new functionals which will be preserved (or
decreased) under the flow with a uniform lower bound, then using
the principle of concentration  to  attack the
compactness/convergence problem. Following our work
\cite{chen994}, M. Struwe \cite{Struwe2000} gave a more concise
proof on Ricci flow and Calabi flow in Riemann surface. This
simple idea of using integral estimates in the flow should be able
to be applied in other geometric flows.

\subsubsection{The Calabi flow on Riemann surfaces}

\vskip-5mm \hspace{5mm}

The Calabi flow is the gradient flow of the K energy and it is a
4th order parabolic equation, proposed by E. Calabi in 1982.
Namely, for a given K\"ahler manifold $(M,[\omega]), $  the Calabi
flow was defined by \[ {{\p \varphi(t)}\over {\p t}} =
 R(\omega_\varphi) - {1\over{vol(M)}} \displaystyle \int_M \; R(\omega) \; \omega^n.
\]
The local existence for this flow is known, while very little is
known for its long term existence since this is a 4th order flow.
The only known result is in Riemann surface where Chrusciel proved
that the flow converges exponentially fast to a unique constant
scalar curvature metric. In 1998 \cite{chen994}, we gave a new
proof based on some geometrical integral estimate and
concentration compactness principle.  Now the challenging
question is:
\begin{q} \rm Does the Calabi flow exists globally for any smooth
initial metric?
\end{q}


\subsubsection{The K\"ahler Ricci flow}

\vskip-5mm \hspace{5mm}

A K\"ahler Ricci flowis defined by
\[ {{\p}\over {\p t}} \omega_\varphi = \omega_\varphi -
Ric(\omega_\varphi).
\]
This flow was first studied by H. D. Cao , following the work of
R. Hamilton on the Ricci flow\footnote{The Ricci flow was
introduced by R. Hamilton \cite{Hamilton82} in 1982. There are
extensive study in this subject (cf. \cite{Hamilton93}) since his
famous work in 3-dimensonal manifold with positive Ricci curvarure
(cf. \cite{Hamilton93} for further references). Another important
geometric flow is the so called ``mean curvature flow.\;" The
codimension 1 case was studied extensively by G. Huisken and many
others. Recently, there are some interesting progress made in
codimension 2 case (cf. \cite{cl1}
 \cite{wang01} for further references).
 }.  Cao\cite{Cao85} proved that the flow always
exists for all the time along with some other interesting results.
It was proved by S. Bando \cite{Bando84} for 3-dimensional
K\"ahler manifolds and by N. Mok \cite{Mok88} for higher
dimensional K\"ahler manifolds that the positivity of bisectional
curvature is preserved under the K\"ahler Ricci flow. The main
issue here is the global convergence on manifold with positive
bisectional curvature. In the work with Tian, we found a set of
new functionals $\{E_k\}_{k=0}^n $ on curvature tensors such that
the Ricci flow is the gradient like flow of these functionals. On
K\"ahler-Einstein manifold with positive scalar curvature, if the
initial metric has positive bisectional curvature, we can prove
that these functionals have a uniform lower bound, via the
effective use of Tian's inequality. Consequently, we are able to
prove the following theorem:

\begin{theo}{\rm \cite{chentian001},\cite{chentian002}}  Let $M$ be a K\"ahler-Einstein
manifold with positive scalar curvature. If the initial metric has
nonnegative bisectional curvature and positive at least at one
point, then the K\"ahler Ricci flow will converge exponentially
fast to a K\"ahler-Einstein metric with constant bisectional
curvature.
\end{theo}
The above theorem in complex dimension $1$ was proved first by
Hamilton \cite{Hamilton88}.   B. Chow \cite{Chow91} later showed
that the assumption that the initial metric has positive curvature
in $S^2$ can be removed since the scalar curvature will become
positive after finite time anyway.

\begin{cor}
The space of K\"ahler metrics with non-negative bisectional
curvature is path-connected.

\end{cor}

Moreover, we can carry over the proof of Theorem 0.12 to a more
general case of K\"ahler orbifolds, for which we will not go into
details here.  Now the definition of these functionals $E_k =
E_k^0 - J_k (k=0,1,\cdots, n)$:

\begin{defi} For any $k=0,1,\cdots, n$, we define a functional $E_k^0$
on ${\cal H}$ by
\[
  E_{k,\omega}^0 (\varphi) = {1\over vol(M)}\; \displaystyle \int_M\;  \left( \log {{\omega_{\varphi}}^n \over \omega^n}
   - h_{\omega}\right) \left(\displaystyle \sum_{i=0}^k\; {{\rm Ric}(\omega_{\varphi})}^{i}\wedge\omega^{k-i} \right)
    \wedge {\omega_{\varphi}}^{n-k} + c_k,
\]
where
\[
c_k ={1\over vol(M)}\; \displaystyle \int_M\;
   h_{\omega} \left(\displaystyle \sum_{i=0}^k\; {{\rm Ric}(\omega)}^{i}\wedge\omega^{k-i} \right)
    \wedge {\omega}^{n-k},
\]
and
\[
   Ric(\omega) -\omega = \sqrt{-1} \p \bar \p h_{\omega},\qquad
   {\rm and}\; \displaystyle\; \int_M (e^{h_{\omega}} - 1 ) \omega^n = 0.
\]
\end{defi}

\begin{defi}
For each $k=0,1,2,\cdots, n-1$, we will define $J_{k, \omega}$ as
follows: Let $\varphi(t) $ ($t\in [0,1]$) be a path from $0$ to
$\varphi$ in ${\cal H}$, we define
\[ J_{k,\o}(\varphi) = -{n-k\over vol(M)}
\int_0^1 \int_M {{\partial \varphi}\over{\partial t}}
\left({\omega_{\varphi}}^{k+1} - {\omega}^{k+1}\right)\wedge
{\omega_{\varphi}}^{n-k-1}\wedge dt.
\]
Put $J_n =0 $ for convenience in notations.
\end{defi}

Note that $E_0$ is the well known K energy function introduced by
T. Mabuchi in 1987.   Direct computations lead to
\begin{theo}For any $k=0,1,\cdots, n$, we have
\begin{eqnarray}
{d E_k \over dt} & = & {{k+1}\over vol(M)} \displaystyle \int_M
\Delta_{\varphi}\left( {{\partial \varphi}\over {\partial t}}
\right )\; {{\rm Ric}(\omega_{\varphi})}^{k} \wedge
{\omega_{\varphi}}^{n-k}
\nonumber \\
& & \qquad -
 {{n-k}\over vol(M)}\displaystyle \int_M {{\partial \varphi}\over {\partial t}} \left({{\rm Ric}(\omega_{\varphi})}^{k+1}
 - {\omega_{\varphi}}^{k+1}\right) \wedge  {\omega_{\varphi}}^{n-k-1}.
\label{eq:decay functional0}
\end{eqnarray}
Here $\{\varphi(t)\}$ is any path in ${\cal H}$.
\end{theo}

Note that under the K\"ahler Ricci flow, these functionals
essentially decreases!  We then prove the derivative of these
functionals along a curve of holomorphic automorphisms give rise
to a set of holomorphic invariants  $\Im_k (k=0,1,\cdots, n)$ (cf.
Theorem 0.1). In case of K\"ahler-Einstein manifolds, all these
invariants vanishes. This give us freedom to re-adjust the flow so
that the evolving K\"ahler potentials are perpendicular to the
first eigenspace of a fixed K\"ahler-Einstein metric. Then we will
be able to show that the evolved volume form has a uniform lower
bound. From this point on, the boot-strapping process will give us
necessary estimates to obtain global convergence.

\vspace*{-0.5mm}

\subsection{Some new result with G. Tian}

\vskip-5mm \hspace{5mm}

In 2001,  Donaldson  proved the following
\begin{theo}{\rm \cite{dona01}} (Openness) For any smooth solution to the geodesic
equation with a disc domain, there are always exists a smooth
solution to the geodesic equation if we perturb the boundary data
in a small open set (of the given boundary data).
\end{theo}
This is somewhat surprising result since it is very hard to deform
any solution of a homogenous Monge-Ampere equation even locally.
However, Donaldson was able to make clever use of the Fredholm
theory of holomorphic discs with totally real boundary in his
proof. Then the problem of closed-ness  becomes very important in
light of this theorem. Tian and I are able to establish the
closed-ness in this case.
\begin{theo}{\rm  \cite{chentian003}} (Closure property) The defomation of geodesic solution
in the preceding theorem is indeed closed, provided we allow
solution to be smooth almost everywhere.
\end{theo}
This is a deep theorem and we will not go into detail here due to
the expository nature of this talk.  However, this theorem, along
with the ideas of proof, shall have implication in both
geometry and other Monge-Ampere type equation in the future.\\

\nocite{cl2} \nocite{clt}

\label{lastpage}

\end{document}